\documentclass[11pt]{article}
\usepackage{caption2}
\usepackage{pdfsync}
\usepackage{amssymb}
\usepackage{color}
\usepackage{amsmath,amsthm}
\usepackage{graphicx}
\usepackage{empheq}
\usepackage{indentfirst}
\usepackage{hyperref}
\newtheorem{definition}{Definition}[section]

\newtheorem{remark}{Remark}[section]
\newtheorem{lem}{Lemma}[section]

\def\f{\frac}
\def\pa{\partial}

\def\e{\eqref}

\def\vep{\varepsilon}

\def\i1n{i=1,\cdots,n}
\def\j1n{j=1,\cdots,n}
\def\ij1n{i,j=1,\cdots,n}

\newcommand{\be}{\begin{equation}}
\newcommand{\ee}{\end{equation}}
\newcommand{\beq}{\begin{equation*}}
\newcommand{\eeq}{\end{equation*}}
 \numberwithin{equation}{section}

\newtheorem{thm}{Theorem}[section]
\DeclareMathOperator*{\esssup}{ess\,sup}
\title{Well-posedness and exact controllability of  the mass balance equations for an extrusion process}
\author{
Mamadou DIAGNE\thanks{Department of Mechanical and Aerospace Engineering,
University of California, San Diego,
La Jolla, CA, 92093, USA.
 E-mail: {\tt mdiagne@eng.ucsd.edu}}
\ Peipei SHANG\thanks{Department of Mathematics, Tongji University, Shanghai 200092, China.
E-mail: {\tt peipeishang@hotmail.com}.}
\ and Zhiqiang WANG\thanks{School of Mathematical Sciences, Fudan
University, Shanghai 200433, China. E-mail: {\tt wzq@fudan.edu.cn}. } 
}
\date{}
\begin{document}

\maketitle

\begin{abstract}
In this paper, we study the well-posedness and exact controllability of a physical model for
a food extrusion process in the isothermal case.
The model expresses the mass balance in the extruder chamber and consists of a hyperbolic Partial Differential Equation (PDE)
 and a nonlinear Ordinary Differential Equation (ODE) whose
 dynamics describes the evolution of a moving interface.
By suitable change of coordinates and fixed point arguments,
we prove the existence, uniqueness and regularity of the solution, and finally the exact controllability  of the coupled system.
\\
\end{abstract}
extruder model.\\
{\bf 2010 MR Subject Classification:}\quad
         35L65, 
         35Q79,  
         35R37,  
         95B05 

\section{Introduction}

The analysis  of
moving boundary problems has been an active subject in the last
decades and their mathematical understanding continues to be an
important interdisciplinary topic  for various engineering  applications.
Representive complex physical systems  describing  biological phenomena  and
reaction diffusion processes such as  Stefan problems  in
cristal growth processes are still calling for many open questions related to their exact  controllability and the design of highly efficient  output feedback control laws for stabilization purposes.  Among these challenging  problems, one can mention the  swelling nanocapsules studied in  \cite{Bouchemal06}, the lyophilization process applied to  pharmaceutical industry
\cite{Daraoui2010,Velardi08}, the cooking processes describing the volume change in food material \cite{Purlis10},
the mixing systems (model of torus reactor including a well-mixed zone and a
transport zone) and the Diesel Oxidation Catalyst (DOC) presented in \cite{Petit10}.

In general, the key point in   such problems is to find    a suitable
change of coordinates which
transforms the system with moving interface  into a system
defined on  a fixed domain. While leading to more complicated infinite dimensional operator,  such approaches allow to utilize the existing results that
establish the existence, uniqueness, regularity and continuous
dependence of solutions  for many   systems of conservation laws.  For the well-posedness problems, we refer to \cite{Stefano05,BressanBook,Godlewski_ESAIM05,Godlewski04,LeFlochBook,LiuYang} (and the
references therein) in the content of weak solutions to systems
(including scalar case) of conservation laws, and
\cite{LiBook94,LiYuBook} in the content of classical solutions to
general quasi-linear hyperbolic systems.
The control problems for hyperbolic conservation laws have been widely studied for a long time.
For controllability of linear hyperbolic systems, one can see the important survey
\cite{Russell}.
The controllability of nonlinear hyperbolic equations (or systems) are studied in \cite{coron, CGWang, Gugat, LiBook09, LiRao}.
Moreover, \cite{CoronBook} provides a comprehensive survey of controllability and stabilization in PDEs
that also includes nonlinear conservation laws.

In this paper we  consider the well-posedness and exact controllability of the Cauchy problem for
a physical model of the extrusion process which  describes
  the mass transport phenomena  in  an isothermal  extruder chamber. Mathematically, the process is described by a hyperbolic PDE defined on a time-varying domain. The dynamics  of the spatial domain  is governed by an
ODE expressing the conservation of mass
in the  extruder and the PDE expresses the convection phenomenon due to the rotating screw.
More detailed description of the model is given in
Section \ref{description-model}.
We mention that  the first result concerning the
mathematical analysis  of the extrusion model as transport equations
coupled via complementary time varying domains is proposed in
\cite{DiagneJESA11}, where the well-posedness for the linearized
model of the extruder is obtained by using perturbation
theory on the linear operator.

Our proof of the well-posedness  (see {\bf Theorem  \ref{thm-well}}) of the  Cauchy
problem for the extrusion model relies on
the characteristic method and the  fixed point  arguments.
The $H^2$-regularity of the solution is proved as well (see {\bf Theorem
\ref{thm-regu}}), which is useful when one considers the asymptotic
stabilization of the corresponding closed-loop system with
feedback controls \cite{DSW2015}.
In this context the   stabilization of a non-isothermal food extrusion process  including the temperature and the moisture content dynamics is investigated in \cite{DSW2015b}.
The idea to prove {\bf Theorem \ref{thm3}} is to construct a solution to \eqref{eq-filling-ratio-bound}-\eqref{eql-normal} which also satisfies the final conditions.
 The way of such construction is based on the controllability result of the linearized system together with fixed point arguments (see for example \cite{CW2012}).

The organization of this paper is as follows: First in Section
\ref{description-model}, we give a description of the extrusion
process model  which is derived from conservation laws.
The main results
 ({\bf Theorems \ref{thm-well}, \ref{thm-regu}, \ref{thm3}}) concerning
the well-posedness, regularity and the controllability of the normalized system are
presented in Section \ref{main-results},
while their proofs are given in Sections  \ref{sec-proof 1}, \ref{sec-proof 2}, \ref{controllability} respectively.

\section{Description of the extrusion process model} \label{description-model}

\begin{figure}[h]
\begin{center}
\includegraphics[scale=0.8]{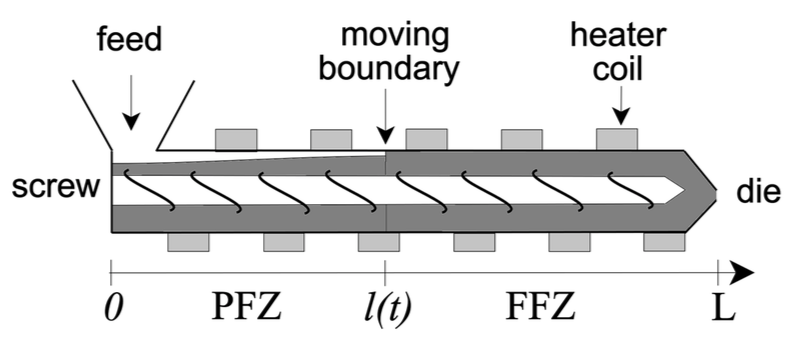}\vspace{0.005cm}
\caption{Schematic description of an extruder} \label{fig1}
\end{center}
\end{figure}

Extruders are designed to process highly viscous materials. They
are mainly used in the chemical industries for polymer processing
as well as in the food industries. An extruder is made of a
barrel, the temperature of which is regulated. One or two
Archimedean screws are rotating inside the barrel. The extruder is
equipped with a die where the material comes out of the process
(see Fig. \ref{fig1}).

In an extruder, the net flow at the die exit is mainly due to the
flow of the material in the  screw axis  direction.
The die resistance  influences highly  the  transport along the extruder and induces an accumulation phenomenon  towards  the down barrel direction.  The accumulation of the material permits us to represent the mass balance in the extruder  using  the \emph{filling ratio} along its  screw channel. More precisely,  the spatial domain of the extruder might be partitioned into a
 \emph{Fully
Filled Zone} ($FFZ$) and a   \emph{Partially
Filled Zone} ($PFZ$). The flow in the \emph{FFZ}
depends on the pressure gradient which appears in this region due to the die
restriction. The $PFZ$ corresponds to   a conveying region that is  submitted to the constant atmospheric pressure and the transport
velocity of the material depends only on the screw speed and pitch. These
two zones are coupled by an interface which is located at the
spatial coordinate  where the pressure gradient changes from zero to a nonzero value. Basically, the
moving interface evolves as a function of
the difference between the feed and die rates.
In the sequel, the spatial domain of the extruder will be taken as
the real interval $\left[0,\, L\right]$, where $L>0$ is the length
of the extruder. Let us denote by $l(t)\in\left[0,\, L\right]$ the
position of the thin interface, the domain of the \emph{PFZ}
is then $\left[0,\, l(t)\right]$ and the \emph{FFZ} is defined on $\left[l(t),\, L\right]$ (see Fig. \ref{fig1}).

Considering the following  change of variables \cite{Diagne}
\be \label{change-var}
x \mapsto y =\f{x}{l(t)}\  \text{in \emph{PFZ}}
\  \text{and} \
x\mapsto y=\f{x-l(t)}{L-l(t)} \  \text{in \emph{FFZ}}
\ee
respectively, the time varying domains $[0,l(t)]$ and $[l(t),L]$ can be transformed to the fixed domain $[0,1]$ in space.
For the sake of simplicity, we still denote by $x$ the space variable instead of $y$.
More precisely, we consider the problem for the corresponding normalized
system defined on $(0,T)\times (0,1)$.

Defining the filling ratio along the $PFZ$  spatial domain, namely, $f_p(t,x)$ as a dynamical variable  \cite{KULSH92,CHIN01},
then the mass balance in this region is  written as follows:
\be\label{eq-filling-ratio-bound}
\begin{cases}
\partial_t f_p(t,x)+\alpha_p \partial_x f_p(t,x)=0, \quad &\text{in}\ \mathbb{R}^+\times(0,1), \\
f_p(0,x)=f^{0}_p(x),\quad &\text{in}\ (0,1),\\
f_p(t,0)=\frac{F_{in}(t)}{\rho_{o} V_{eff}N(t)},\quad &\text{in}\ \mathbb{R}^+,
\end{cases}
\ee
where
\begin{align}\label{alphap-normal}
\alpha_p(x,N(t),l(t),f_p(t,1))=\f{\zeta N(t)-x F(l(t),N(t),f_p(t,1))}{l(t)}
\end{align}
and $\alpha_p$ is the transport velocity of the material, $\zeta$ is the screw pitch,
$F_{in}$ denotes the feed rate, $\rho_0$ is the melt density,
$V_{eff}$ is the effective volume and $N(t)$ is the rotation speed of the screw.
$F$ is the dynamics of the moving interface described by equations \eqref{eql-normal}-\eqref{F}.
 The interface motion is generated by the gradient
 of pressure which appears in the \emph{FFZ} \cite{KIM001,KIM002,KULSH92,CHIN01} and
 under the assumption of constant viscosity along the extruder (the isothermal case), its evolution is given by the following
\be\label{eql-normal}
\begin{cases}
\dot{l}(t)=F(l(t),N(t),f_p(t,1)),\quad &\text{in}\ \mathbb{R}^+,\\
l(0)=l^0,
\end{cases}
\ee
where
\be\label{F}F(l(t),N(t),f_p(t,1))=N(t)g(l(t),f_p(t,1)),
\ee
with
\be\label{gcon}
g(l(t),f_p(t,1))=\frac{\zeta K_d\left(L-l(t)\right)}{\left[B\rho_0+K_d (L-l(t))\right](1-f_p(t,1))}-\frac{ \zeta f_p(t,1)}{1-f_p(t,1)}.
\ee
In \eqref{gcon}, $K_d$
denotes the die conductance and $B$ is the geometric parameter.

In the whole paper, unless otherwise specified,
we always assume that $l^0\in(0,L)$, $f^0_p\in W^{1,\infty}(0,1)$, $F_{in},N\in L^{\infty}(0,T)$,
$F_{in}/N\in W^{1,\infty}(0,T)$.
For the sake of simplicity, we denote from now on
$\|f\|_{L^{\infty}}$ ($\|f\|_{W^{1,\infty}}$, $\|f\|_{L^2}$, resp.) as the $L^{\infty}$ ($W^{1,\infty}$, $L^2$, resp.)
norm of the function $f$ with respect to its variables.

\section{Main Results}\label{main-results}

In this section, we present the main results on the well-posedness and exact controllability of
the coupled system \eqref{eq-filling-ratio-bound}-\eqref{eql-normal},
we have the following two theorems.

\begin{thm}\label{thm-well}
Let $T>0$ and $(l_e, N_e, f_{pe})$ be a constant equilibrium, i.e.,
\be\label{equil}
F(l_e,N_e,f_{pe})=0
\ee
with $0<f_{pe}<1$, $0<l_e<L$.
Assume that the compatibility condition at $(0,0)$ holds
\be \label{compati}
\f{F_{in}(0)}{\rho_{0}V_{eff}N(0)}=f^0_p(0).
\ee
Then, there exists $\varepsilon_0>0$ (depending on $T$) such that
for any $\varepsilon\in(0,\varepsilon_0]$,
if
\begin{eqnarray}
|l^0\!-\!l_e|+\|f^{0}_p(\cdot)-f_{pe}\|_{W^{1,\infty}}
+ \|\frac{F_{in}(\cdot)}{\rho_{0}V_{eff}N(\cdot)}- f_{pe}\|_{W^{1,\infty}}
+ \|N(\cdot)\!-\!N_e\|_{L^{\infty}}\!\leq \!\varepsilon, \label{G3}
\end{eqnarray}
Cauchy problem \eqref{eq-filling-ratio-bound}-\eqref{eql-normal}
admits a unique solution
$(l,f_p)\!\in \!W^{1,\infty}(0,T)\!\times \!W^{1,\infty}((0,T)\!\times\!(0,1))$,
and the following estimates hold
\begin{align} \label{est-fp}
\|l(\cdot)-l_e\|_{W^{1,\infty}}&\leq  C_{\varepsilon_0}\cdot\varepsilon,\\
\|f_p(\cdot,\cdot)-f_{pe}\|_{W^{1,\infty}}&\leq  C_{\varepsilon_0}\cdot\varepsilon,
\end{align}
where $C_{\varepsilon_0}$ is a constant depending on $\varepsilon_0$,
but independent of $\varepsilon$.
\end{thm}

\begin{thm}\label{thm-regu}
Under the assumptions of Theorem \ref{thm-well}, we assume furthermore that
$f^0_p(\cdot)\in H^2(0,1)$, $\displaystyle\f{F_{in}(\cdot)}{\rho_{0}V_{eff}N(\cdot)}\in H^2(0,T)$,
and the compatibility condition at $(0,0)$ holds
\be \label{compati2}
(f^0_{p})_x(0)+\displaystyle\f{l(0)}{\zeta N(0)}\cdot \f{d}{dt} \Big( \f{F_{in}(t)}{\rho_{0}V_{eff}N(t)} \Big)\Big |_{t=0}=0.
\ee
Then, there exists $\varepsilon_0>0$ (depending on $T$) such that for any $\varepsilon\in (0,\varepsilon_0]$, if
\begin{eqnarray}
 |l^0\!-\!l_e|+\|f^{0}_p(\cdot)\!-\!f_{pe}\|_{H^2(0,1)}
  + \|\frac{F_{in}(\cdot)}{\rho_{0}V_{eff}N(\cdot)} - f_{pe}\|_{H^2(0,T)}
  +\|N(\cdot)\!-\!N_e\|_{L^{\infty}}\!\leq \!\varepsilon,   \label{G333}
\end{eqnarray}
Cauchy problem \eqref{eq-filling-ratio-bound}-\eqref{eql-normal} has a unique solution
$(l,f_p)\in W^{1,\infty}(0,T)\times C^0([0,T];H^2(0,1))$
with the additional estimate
\be
\|f_p(\cdot,\cdot)-f_{p_e}\|_{C^0([0,T];H^2(0,1))}\leq  C_{\varepsilon_0}\cdot \varepsilon,
\ee
where $C_{\varepsilon_0}$ is a constant depending on $\varepsilon_0$,
but independent of $\varepsilon$.
\end{thm}

\begin{remark}
The solution in Theorem \ref{thm-well} or in Theorem \ref{thm-regu}
is often called semi-global solution since it exists on any preassigned time interval
$[0,T]$ provided that $(l,f_p)$ has some kind of smallness (depending on $T$), see \cite{Li2001, W2006}.
\end{remark}

\begin{remark}\label{rem}
We have the hidden regularity that $f_p \in C^0([0,1];H^2(0,T))$ in Theorem \ref{thm-regu} (see \cite{CKW, SW} for the idea of proof).
\end{remark}

The problem of exact  controllability for Cauchy problem \eqref{eq-filling-ratio-bound}-\eqref{eql-normal} can be described as follows:
For any given initial data $(l^0, f^0_p(x))$, any final data
$(l^1,f^1_p(x))$, to find a time $T$ and controls $F_{in}(t)$ and $N(t)$
such that the solution to Cauchy problem \eqref{eq-filling-ratio-bound}-\eqref{eql-normal} satisfies
\begin{align}\label{final-con1}
l(T)&=l^1,\\
\label{final-con2}
f_p(T,x)&=f^1_p(x).
\end{align}
Our result is  the following theorem on local controllability  in the sense that the initial and final data are both close to the given equilibrium determined by
\eqref{equil}.

\begin{thm} \label{thm3}
Let
\be\label{Te}
T_e:=\f{l_e}{\zeta N_e}
\ee
be the critical control time. Then, for any $T>T_e$,
there exists $\nu_1>0$ suitably small such that, for any $\nu\in(0,\nu_1]$, $l^0,l^1\in (0,L)$ and $f^0_p, f^1_p\in W^{1,\infty}(0,1)$ with
\begin{align} \label{assum3}
&|l^0\!-\!l_e|+|l^1-l_e|+\|f^{0}_p(\cdot)-f_{pe}\|_{W^{1,\infty}}+\|f^1_p(\cdot)-f_{pe}\|_{W^{1,\infty}}
\leq \nu,
\end{align}
there exist $N\in L^{\infty}(0,T)$ and $F_{in} \in L^{\infty}(0,T)$ satisfying
\be\label{nFN}
\|\frac{F_{in}(\cdot)}{\rho_{0}V_{eff}N(\cdot)}- f_{pe}\|_{W^{1,\infty}}
+ \|N(\cdot)\!-\!N_e\|_{L^{\infty}}
\leq C_{\nu_1}\cdot\nu,
\ee
such that the weak solution $(l(t),f_p(t,x))$ to Cauchy problem \eqref{eq-filling-ratio-bound}-\eqref{eql-normal} satisfies the final condition
\eqref{final-con1}-\eqref{final-con2}. Here, $C_{\nu_1}$ is a constant depending on $\nu_1$, but independent of $\nu$.
\end{thm}


\section{Proof of Theorem \ref{thm-well}} \label{sec-proof 1}

In order to conclude Theorem \ref{thm-well}, it suffices to prove the following lemma on
local well-posedness of  Cauchy problem \eqref{eq-filling-ratio-bound}-\eqref{eql-normal}.

\begin{lem}\label{lem-loc}
There exist $\varepsilon_1>0$ and $\delta>0$ suitably small, such that for any $\varepsilon\in(0,\varepsilon_1]$,
$l^0\in (0,L)$, $f^0_p\in W^{1,\infty}(0,1)$, $F_{in}/ N\in W^{1,\infty}(0,T)$ with
\be\label{initial-bound}
\!|l^0\!-\!l_e|+\|f^{0}_p(\cdot)\!-\!f_{pe}\|_{W^{1,\infty}}\!
   +\!\|\frac{F_{in}(\cdot)}{\rho_{0}V_{eff}N(\cdot)}\!-\!f_{pe}\|_{W^{1,\infty}}\!
   +\!\|N(\cdot)\!-\!N_e\|_{L^{\infty}}\!\!\leq \!\varepsilon,
\ee
Cauchy problem \eqref{eq-filling-ratio-bound}-\eqref{eql-normal}
admits a unique local solution on $[0,\delta]$, which satisfies
the following estimates
\begin{align}
\label{estimatel}
|l(t)-l_e|&\leq  C_{\varepsilon_1}\cdot\varepsilon,\quad \forall t\in [0,\delta],\\
\label{estimatefp}
\|f_p(t,\cdot)-f_{pe}\|_{W^{1,\infty}}
&\leq  C_{\varepsilon_1}\cdot\varepsilon,\quad \forall t\in[0,\delta],
\end{align}
where $C_{\varepsilon_1}$ is a constant depending on $\varepsilon_1$, but independent of $\varepsilon$.
\end{lem}

Let us first show how to conclude Theorem \ref{thm-well} from Lemma \ref{lem-loc}.
By Lemma \ref{lem-loc}, we take $\varepsilon_2 \in (0, \varepsilon_1]$
such that $C_{\varepsilon_1}\cdot\varepsilon_2\leq  \varepsilon_1$.
Then  for any $\varepsilon\in(0,\varepsilon_2]$ and any initial-boundary data such that
\eqref{initial-bound} holds,
Cauchy problem \eqref{eq-filling-ratio-bound}-\eqref{eql-normal} admits a unique solution on $[0,\delta]$.
Furthermore, one has
\begin{align}
\label{n-use2}
|l(\delta)-l_e|&\leq  C_{\varepsilon_1}\cdot \varepsilon\leq  \varepsilon_1,\\
\label{n-use}
\|f_p(\delta,\cdot)-f_{pe}\|_{W^{1,\infty}}&\leq  C_{\varepsilon_1}\cdot\varepsilon\leq  \varepsilon_1.
\end{align}
By taking $(l(\delta),f_p(\delta,\cdot))$ as new initial data  and applying Lemma \ref{lem-loc} on $[\delta,2\delta]$,
the solution of Cauchy problem \eqref{eq-filling-ratio-bound}-\eqref{eql-normal} is extended to $[0,2\delta]$.
For fixed $T>0$, we can extend the local solution to Cauchy problem
 \eqref{eq-filling-ratio-bound}-\eqref{eql-normal} to $[0,T]$ eventually
by reducing the value of $\varepsilon_0$ and applying Lemma \ref{lem-loc}
in finite times (at most $[T/\delta]+1$ times).
Therefore, to conclude Theorem \ref{thm-well}, it remains to prove Lemma  \ref{lem-loc}.  \qed


\noindent
{\bf Proof of Lemma \ref{lem-loc}:}
The proof is divided into 4 steps.

\noindent \textbf{Step 1. Existence and uniqueness of $(l(\cdot),f_p(\cdot,1))$ by fixed point argument.}

Let $\varepsilon_1>0$ be such that
\be\label{varepsilon1}
0<\varepsilon_1<\min\{l_e,L-l_e,f_{pe},1-f_{pe}\}.
\ee
Denote
\begin{align}
\label{FW}
\|F\|_{W^{1,\infty}}:&=\sum_{|\alpha|\leq  1}\esssup_{\substack{0<x_1<L\\
N_e-\varepsilon_1<x_2<N_e+\varepsilon_1\\
0<x_3<1}}|D^{\alpha}F(x_1,x_2,x_3)|,\\
\Psi(t):&=(l(t),f_p(t,1)),\quad t\in[0,T].
\end{align}

\begin{figure}[htbp!]
\includegraphics[width=0.5\textwidth]{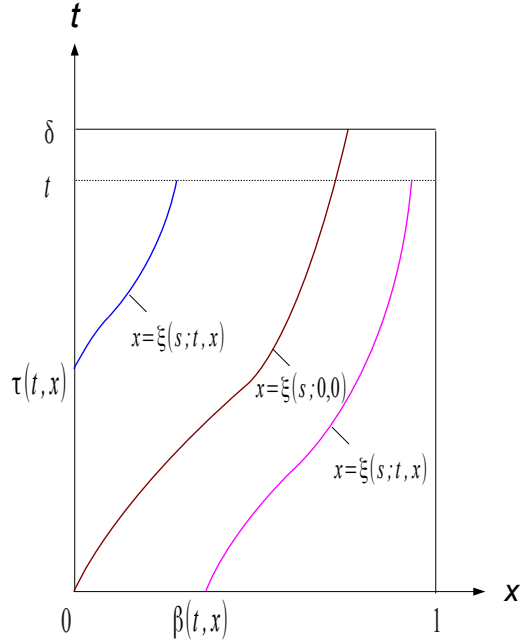}
\centering
\caption{The characteristics $\xi(s;t,x)$ and $\beta(t,x)$, $\tau(t,x)$. \label{fig2} }\par\vspace{0pt}
\end{figure}

For any given $\delta>0$ small enough (to be chosen later), we define
a domain candidate as a closed subset of $C^0([0,\delta])$ with respect to $C^0$ norm:
\begin{align}
\Omega_{\delta,\varepsilon_1}:=
\Big\{ \Psi \in C^0 ([0 ,\delta]) | \Psi(0)=(l^0,f^0_{p}(1)),\
\|\Psi(\cdot)-(l_e,f_{pe})\|_{C^0([0,\delta])} \leq \vep_1 \Big\}.
\end{align}
We denote by $\xi(s;t,x)$, with $(s, \xi(s;t,x)) \in [0,t]\times [0,1]$,
the characteristic curve passing through the point $(t,x) \in [0,\delta] \times [0,1]$ (see Fig. \ref{fig2}), i.e.,
\be\label{odexi2}
\begin{cases}
\displaystyle\frac{d \xi(s;t,x)} {ds}=\alpha_p\bigl(\xi(s;t,x),N(s),l(s),f_p(s,1)\bigr),\\
\xi(t;t,x)=x.
\end{cases}
\ee

Let us define a map $\mathfrak{F}:=(\mathfrak{F}_1,\mathfrak{F}_2)$, where
$\mathfrak{F} : \Omega_{\delta,\vep_1} \rightarrow  C^0 ([0 ,\delta])$, $\Psi\mapsto \mathfrak{F}(\Psi)$
as
\begin{align}\label{map1}
\mathfrak{F}_1(\Psi)(t)&:=l^0+\int_0^t F\bigl(l(s),N(s),f_p(s,1)\bigr)ds,\\
\label{map2}
\mathfrak{F}_2(\Psi)(t)&:=f^0_p(\xi(0;t,1)).
\end{align}
Solving the linear ODE \eqref{odexi2} with $\alpha_p$ given by \eqref{alphap-normal},
one easily gets for all $\delta $ small and all $0\leq  s\leq  t \leq \delta$ that
\be\label{solodexi2}
\xi(s;t,1)=e^{\int_s^t \f{F(l(\sigma),N(\sigma),f_p(\sigma,1))}{l(\sigma)}\, d\sigma}-\int_s^t
\f{\zeta N(\sigma)}{l(\sigma)}\cdot e^{\int^{\sigma}_s \f{F(l(s),N(s),f_p(s,1))}{l(s)}\, ds}\,d\sigma.
\ee
It is obvious that $\mathfrak{F}$ maps into $\Omega_{\delta,\vep_1}$ itself if
\be \label{delta01}
0<\delta< \min \Big\{T,\frac{l_e-\varepsilon_1}{\zeta(N_e+\varepsilon_1)},\f{l_e-\vep_1}{\|F\|_{W^{1,\infty}}},\f{L-l_e-\vep_1}{\|F\|_{W^{1,\infty}}}\Big \},
\ee
%

Now we prove that, if $\delta$ is small enough, $\mathfrak{F}$ is a contraction mapping on
$\Omega_{\delta,\vep_1}$ with respect to the $C^0$ norm.
Let $\Psi=(l,f_p),\bar \Psi=(\bar l,\bar f_p)\in \Omega_{\delta,\vep_1}$.
We denote by $\bar \xi(s;t,x)$ the corresponding characteristic curve  passing through $(t,x)$:
   \be
   \begin{cases}
\displaystyle\frac{d \bar \xi(s;t,x)} {ds}=\alpha_p\bigl(\bar \xi(s;t,x),N(s),\bar l(s), \bar f_p(s,1)\bigr),\\
\bar\xi(t;t,x)=x.
\end{cases}
   \ee
Similarly as \e{solodexi2},  one has for all $\delta $ small and all $0\leq  s\leq  t \leq \delta$ that
\be\label{solodebarxi2}
\bar\xi(s;t,1)=e^{\int_s^t \f{F(\bar l(\sigma),N(\sigma),\bar f_p(\sigma,1))}{\bar l(\sigma)}\, d\sigma}-\int_s^t
\f{\zeta N(\sigma)}{\bar l(\sigma)}\cdot e^{ \int^{\sigma}_s \f{F(\bar l(s),N(s),\bar f_p(s,1))}{\bar l(s)}\,
ds}\,d\sigma.
\ee

Therefore it holds for all $t\in[0,\delta]$ that
 \begin{align}
|\mathfrak{F}_1(\bar \Psi)(t)-\mathfrak{F}_1(\Psi)(t)|=&
\Big|\int_0^t F\bigl(\bar l(s),N(s),\bar f_p(s,1))\, ds-\int_0^t F\bigl(l(s),N(s),f_p(s,1))\, ds \Big |\nonumber\\
\leq  &\delta \|F\|_{W^{1,\infty}} \|\bar \Psi-\Psi\|_{C^0([0,\delta])}.   \label{difference1}
   \end{align}
On the other hand, it follows from \eqref{map2}, \eqref{solodexi2} and \eqref{solodebarxi2},   that for  all $t\in [0,\delta]$,
\begin{align}
      &|\mathfrak{F}_2(\bar \Psi)(t)-\mathfrak{F}_2(\Psi)(t)|\nonumber
        \\
    =&|f^0_p(\bar\xi(0;t,1))-f^0_p(\xi(0;t,1))|\nonumber
       \\
    \leq  & \|f^0_{px}\|_{L^{\infty}}\Big(\Big|e^{\int_0^t \f{F(\bar l(\sigma),N(\sigma),\bar f_p(\sigma,1))}{\bar l(\sigma)}\, d\sigma}
                    -e^{\int_0^t \f{F(l(\sigma),N(\sigma),f_p(\sigma,1))}{l(\sigma)}\, d\sigma}\Big|\nonumber
       \\
      & \qquad \qquad + \int_0^t\Big|\f{\zeta N(\sigma)}{\bar l(\sigma)}\cdot
          e^{\int^{\sigma}_0 \f{F(\bar l(s),N(s),\bar f_p(s,1))}{\bar l(s)}\, ds}-\f{\zeta N(\sigma)}{l(\sigma)}\cdot
          e^{\int^{\sigma}_0 \f{F(l(s),N(s),f_p(s,1))}{l(s)}\, ds}\Big|\, d\sigma \Big).\nonumber
\end{align}
By \e{FW} and the fact that $\Psi,\bar \Psi \in \Omega_{\delta,\vep_1}$ and
$l(t)\geq l_e-\vep_1>0, |N(t)|\leq N_e+\vep_1 $, it follows  that for  all $t\in [0,\delta]$,
\be \label{difference22}
      |\mathfrak{F}_2(\bar \Psi)(t)-\mathfrak{F}_2(\Psi)(t)|
\leq C \delta \|f^0_p(\cdot)-f_{pe}\|_{W^{1,\infty}}  \|\bar \Psi-\Psi \|_{C^0([0,\delta])},
\ee
where $C>0$ is a constant
 independent of $(\bar \Psi,\Psi)$.
Finally, combining \eqref{difference1} and \eqref{difference22}, we can choose $\delta$ small enough
such that
\be \label{estimatedelta0}
       \|\mathfrak{F}(\bar\Psi)-\mathfrak{F}(\Psi)\|_{C^0([0,\delta])}
    \leq   \frac{1}{2}  \|\bar \Psi-\Psi \|_{C^0([0,\delta])}.
\ee
Banach fixed point  theorem implies the existence of the unique fixed point
$(l(\cdot),f_p(\cdot,1))$ of the mapping $\mathfrak{F}$: $\Psi = \mathfrak{F}(\Psi)$
in $\Omega_{\delta,\vep_1}$.

\noindent \textbf{Step 2. Construction of a solution by characteristic method.}

With the existence of $(l(\cdot),f_p(\cdot,1))$ and $\delta>0$ by \textbf{Step 1},
we can construct a solution to Cauchy problem \eqref{eq-filling-ratio-bound}-\eqref{eql-normal}.
For every $(t,x)$ in $[0,\delta]\times [0,1]$, we still denote by $\xi(s;t,x)$, with $(s, \xi(s;t,x)) \in [0,t]\times [0,1]$,
the characteristic curve passing through the point $(t,x)$, see \e{odexi2}.
Since the velocity function $\alpha_p$ is positive, the characteristic $\xi(s;t,x)$ intersects the $x$-axis at point $(0,\beta(t,x))$
with $\beta(t,x)=\xi(0;t,x)$ if  $0\leq  \xi(t;0,0) \leq x\leq 1$;
the characteristic $\xi(s;t,x)$ intersects the $t$-axis at point $(\tau(t,x),0)$ with
$\xi(\tau(t,x);t,x)=0$ if  $0\leq x\leq \xi(t;0,0)$ (see Fig. \ref{fig2}).
Moreover, we have (see \cite[Lemma 3.2 and its proof, Page 90-91]{LiYuBook} for a more general situation)
\begin{align}
\label{applem12}
\f{\pa \tau(t,x)}{\pa x}&=\f{-l(\tau(t,x))}{\zeta N(\tau(t,x))}
  \, e^{\int_{\tau(t,x)}^t \f{F(l(\sigma),N(\sigma),f_p(\sigma,1))}{l(\sigma)}\, d\sigma},
   \\    \label{applem22}
\f{\pa \beta(t,x)}{\pa x}
    &= e^{\int_0^t \f{F(l(\sigma),N(\sigma),f_p(\sigma,1))}{l(\sigma)}\, d\sigma}.
\end{align}
We define $f_p$ by
\be\label{sol1}
f_p(t,x)=
\begin{cases}
\displaystyle\frac{F_{in}(\tau(t,x))}{\rho_{0}V_{eff}N(\tau(t,x))},
            \quad \text{if}\ 0\leq  x \leq  \xi(t;0,0) \leq 1, 0\leq t\leq \delta,
             \\
 f^{0}_p(\beta(t,x)),\quad \text{if}\  0\leq  \xi(t;0,0)\leq  x\leq  1, 0\leq t\leq \delta.
 \end{cases}
\ee
Then it is easy to check that
$(l,f_p)\in W^{1,\infty}(0,\delta) \times W^{1,\infty}((0,\delta)\times (0,1))$
under the compatibility condition \e{compati} and $(l,f_p)$ is indeed a solution
to Cauchy problem \eqref{eq-filling-ratio-bound}-\eqref{eql-normal}.

\noindent \textbf{Step 3. Uniqueness of the solution.}

Assume that Cauchy problem \eqref{eq-filling-ratio-bound}-\eqref{eql-normal}
has two solutions $(l,f_p),(\bar l,\bar f_p)$ on $[0,\delta] \times [0,1]$.
 It follows that $(l(\cdot),f_p(\cdot,1))=(\bar l(\cdot),\bar f_p(\cdot,1))$
since they are both the fixed point of the mapping $\mathfrak{F}$: $\Psi = \mathfrak{F}(\Psi)$
in $\Omega_{\delta,\vep_1}$. This fact implies that the characteristics $\xi(\cdot;t,x)$ and
 $\bar \xi(\cdot;t,x)$ coincide with each other
and therefore so do the solutions $f_p$ and $\bar f_p$ by characteristic method.

\noindent \textbf{Step 4. A priori estimate on the local solution.}

By definition of $f_p$ and assumption  \eqref{initial-bound}, it is clear that  for all $t\in[0,\delta]$,
\be
\|f_p(t,\cdot)-f_{pe}\|_{L^{\infty}}\leq  \varepsilon. \label{G4}
\ee
Thanks to \e{eql-normal}, \e{equil}, \e{G4} and assumption \e{initial-bound},
we get for all $t\in [0,\delta]$ that
\begin{align*}
             |\dot{l}(t)|
      =&|F\left(l(t),N(t),f_p(t,1)\right)-F\left(l_e,N_e,f_{pe}\right)|\nonumber
             \\ \nonumber
    \leq & \|F\|_{W^{1,\infty}}  (|l(t)-l_e|+|N(t)-N_e|+|f_p(t,1)-f_{pe}|).
            \\
     \leq  &   \|F\|_{W^{1,\infty}} |l(t)-l_e|+ 2 \vep \|F\|_{W^{1,\infty}},
   \end{align*}
which yields \eqref{estimatel} from \eqref{initial-bound} and Gronwall's inequality.
On the other hand, from \eqref{sol1}
\begin{align}\label{partialfpx}
        \Big\|\f{\pa f_p}{\pa x}\Big\|_{L^{\infty}}
   \leq  &\Big\|\f{\pa }{\pa x} \Big(\frac{F_{in}(\tau(t,x))}{\rho_{0}V_{eff}N(\tau(t,x))}\Big) \Big\|_{L^{\infty}}
                   +\Big\|\f{\pa }{\pa x} f^0_p(\beta(t,x)) \Big\|_{L^{\infty}}\nonumber
           \\
    \leq  &\Big\|\frac{F_{in}(\cdot)}{\rho_{0}V_{eff}N(\cdot)}-f_{pe}\Big\|_{W^{1,\infty}}
                      \Big\|\f{\pa \tau}{\pa x}\Big\|_{L^{\infty}}
      +\Big\|f^0_p(\cdot)-f_{pe}\Big\|_{W^{1,\infty}} \Big\|\f{\pa \beta}{\pa x}\Big\|_{L^{\infty}}.
\end{align}
Combining \eqref{applem12}, \eqref{applem22}, \eqref{partialfpx} and assumption \eqref{initial-bound},
we obtain \eqref{estimatefp} which concludes
the proof of Lemma \ref{lem-loc}.    \qed


\section{Proof of Theorem \ref{thm-regu}}
\label{sec-proof 2}
Before proving Theorem \ref{thm-regu},
let us recall a classical result on Cauchy problem of the following general linear transport equation
\be\label{equation-u}
\begin{cases}
u_t+a(t,x)u_x=b(t,x)u+c(t,x),\quad (t,x)\in (0,T)\times(0,1),\\
u(0,x)=u_0(x),\quad x\in (0,1),\\
u(t,0)=h(t),\quad t\in (0,T),
\end{cases}
\ee
where $a(t,x) > 0$, $a, a_x, b\in L^{\infty}((0,T)\times(0,1))$ and $c\in L^2((0,T)\times(0,1))$.

We recall from \cite[Section 2.1]{CoronBook}, the definition of a weak solution to Cauchy problem
\eqref{equation-u}.
\begin{definition}\label{definition-u}
Let $T>0$, $u_0\in L^2(0,1)$, $h\in L^2(0,T)$ be given.
A weak solution of Cauchy problem
\eqref{equation-u} is a function
$u\in C^0([0,T];L^2(0,1))$
such that for every $\tau\in[0,T]$, every test function $\varphi\in C^1([0,T]\times[0,1])$ such that
$ \varphi(t,1)=0,\ \forall t\in[0,T]$,
one has
\begin{align}\label{definition-u1}
&- \int_0^{\tau}\int_0^1 \Big(u [\pa_t\varphi+a \pa_x \varphi+(a_x+b)\varphi ]  +  c\varphi \Big)\, dx\, dt
   + \int_0^1 u(\tau,\cdot) \varphi(\tau,\cdot)  \, dx   \nonumber\\
&- \int_0^1 u_0\varphi(0,\cdot)\, dx
- \int_0^{\tau}ha(\cdot,0)\varphi(\cdot,0)\, dt=0.
\end{align}
\end{definition}
We have the following lemma
\begin{lem}\label{lem-lin}
Let $T>0$, $u_0 \in L^2(0,1)$ and $h \in L^2(0,T)$ be given.
Then, Cauchy problem \eqref{equation-u} has
a unique weak solution $u \in C^0([0,T];L^2(0,1))$ and the following estimate holds:
\begin{align}\label{ape}
\|u\|_{C^0([0,T];L^2(0,1))}\leq  C (\|u_0\|_{L^2(0,1)}+\|h\|_{L^2(0,T)}+\|c\|_{L^2((0,T)\times(0,1))} ),
\end{align}
where $C=C(T, \|a\|_{L^{\infty}((0,T)\times(0,1))},\|a_x\|_{L^{\infty}((0,T)\times(0,1))},\|b\|_{L^{\infty}((0,T)\times(0,1))})$
is a constant independent of $u_0,h,c$.
\end{lem}

For the proof of Lemma \ref{lem-lin}, one can refer to \cite{LiYuBook} for classical solution or
\cite[Theorem 23.1.2, Page 387]{Horma3} for Cauchy problem on $\mathbb{R}$ without boundary.

\noindent
{\bf Proof of Theorem \ref{thm-regu}.}
By Theorem \ref{thm-well} and  Lemma \ref{lem-lin}, it suffices to prove that  the systems of $f_{p_{xx}}$
satisfies all the assumptions of Lemma \ref{lem-lin}.

Differentiating \eqref{eq-filling-ratio-bound} with respect to $x$ once and twice give us successively  that
\be\label{de-eq-filling-ratio-bound}
\begin{cases}
\partial_t f_{p_x}(t,x)+\alpha_p(t,x) \partial_x f_{p_x}(t,x)=-\alpha_{p_x}(t,x)  f_{p_x}(t,x),\quad &\text{in}\ (0,T)\times(0,1),\\
f_{p_x}(0,x)=f^{0}_{p_x}(x),\quad &\text{in}\ (0,1),\\
f_{p_x}(t,0)=\displaystyle\f{-l(t)}{\zeta N(t)}\cdot \f{d}{dt}\Big(\f{F_{in}(t)}{\rho_{0}V_{eff}N(t)}\Big),\quad &\text{in}\ (0,T),
\end{cases}
\ee
and
\be\label{de-de-eq-filling-ratio-bound}
\begin{cases}
\partial_t f_{p_{xx}}(t,x)+\alpha_p(t,x)  \partial_x f_{p_{xx}}(t,x)
       =-2\alpha_{p_x}(t,x) f_{p_{xx}}(t,x),\quad &\text{in}\ (0,T)\times(0,1),
         \\
f_{p_{xx}}(0,x)=f^{0}_{p_{xx}}(x),\quad &\text{in}\ (0,1),
       \\
f_{p_{xx}}(t,0)=\displaystyle\f{-l(t)}{\zeta N(t)}
       \Big [\f{F(l(t),N(t),f_p(t,1))}{\zeta N(t)} \cdot \f{d}{dt}\Big(\f{F_{in}(t)}{\rho_{0}V_{eff}N(t)}\Big) \\
       \qquad \qquad \qquad \qquad  \qquad \displaystyle
          -\f{d}{dt}\Big(\f{l(t)}{\zeta N(t)} \cdot \f{d}{dt}\Big(\f{F_{in}(t)}{\rho_{0}V_{eff}N(t)}\Big)\Big) \Big ], \quad &\text{in}\ (0,T),
\end{cases}
\ee
with
     \be  \alpha_p(t,x) =\f{\zeta N(t)-xF(l(t),N(t),f_p(t,1))}{l(t)},
     \quad  \alpha_{p_x}(t,x)=\f{-F(l(t),N(t),f_p(t,1))}{l(t)}.
\ee
From the assumptions that $f_p^0\in H^2(0,1)$, $\displaystyle\f{F_{in}(\cdot)}{\rho_{0}V_{eff}N(\cdot)}\in H^2(0,T)$
and the compatibility conditions \e{compati} and \e{compati2},
one easily concludes Theorem \ref{thm-regu}   by applying Lemma \ref{lem-lin}
  to Cauchy problem \eqref{de-de-eq-filling-ratio-bound}.              \qed

\section{Proof of Theorem \ref{thm3}}\label{controllability}

The idea to prove Theorem \ref{thm3} is to construct a solution to Cauchy problem \eqref{eq-filling-ratio-bound}-\eqref{eql-normal} which also satisfies the final conditions.
 The way of such construction is based on the controllability result of the linearized system together with fixed point arguments (see \cite{CW2012}).

For any fixed initial data $(l^0, f_p^0(x))$ and final data $(l^1, f_p^1(x))$ close to the equilibrium $(l_e, f_{pe})$,
we define a domain as a closed subset of $C^0([0,T])$ with respect to
$C^0$ norm:
\begin{align*}
\Omega^{\varepsilon_1,T}:=\Big\{\Phi\in C^0([0,T]) | & \Phi(0)= (l^0, f_p^0(1)),  \Phi(T)= (l^1, f_p^1(1)), \ \|\Phi(\cdot)-(l_e,f_{pe})\|_{C^0([0,T])}\leq \varepsilon_1\Big\},
\end{align*}
where the constant $\varepsilon_1>0$ is determined by \eqref{varepsilon1}.
We study the exact controllability for the linearized system deduced
from equations \eqref{eq-filling-ratio-bound}-\eqref{eql-normal}, replacing  $(l(t), f_p(t,1))$ by $(a(t),b(t))$ in functions $F$ and $\alpha_p$:
\be\label{eql}
\begin{cases}
\dot{l}(t)=F^{a,b}(t),\quad &\text{in} \ (0,T),\\
\partial_t f_p(t,x)+\alpha_p^{a,b} (t,x) \partial_x f_p(t,x)=0, \quad &\text{in}\ (0,T)\times(0,1), \\
l(0)=l^0,\quad f_p(0,x)=f^{0}_p(x),\quad &\text{in}\ (0,1),\\
\displaystyle f_p(t,0)=\frac{F_{in}^{a,b}(t)}{\rho_{o} V_{eff}N^{a,b}(t)},\quad &\text{in}\ (0,T),
\end{cases}
\ee
where
\begin{align}
F^{a,b}(t)& =   N^{a,b}(t) g^{a,b}(t),
\\
g^{a,b}(t) &= \frac{ \zeta K_d\left(L-a(t)\right)}{\left[B\rho_0+K_d (L-a(t))\right](1-b(t))}-\frac{ \zeta b(t)}{1-b(t)},
\\
\label{talpha_p}
\alpha_p^{a,b}   (t,x)
&= \f{\zeta N^{a,b}(t) - xF^{a,b}(t)}{a(t)}.
\end{align}
By assumption $T>T_e$, it is possible to find controls $N^{a,b}(t)$ and $F_{in}^{a,b}(t)$ such that the solution of \eqref{eql} satisfies \eqref{final-con1}-\eqref{final-con2}.

For any $\Phi=(a,b)\in \Omega^{\varepsilon_1,T}$, we choose the control function $N^{a,b}(t)$ as
\be\label{cN}
N^{a,b}(t):=\f{l^1-l^0}{T}\cdot \f{1}{g^{a,b}(t)}, \quad t\in [0,T],
\ee
such that
\be\label{FF}
F^{a,b}(t) =\f{l^1-l^0}{T}, \quad t\in [0,T],
\ee
thus $l(t)$ is a linear function of $t$:
\be\label{l}
l(t)=l^0+\f{l^1-l^0}{T}t, \quad t\in [0,T].
\ee
By assumption  \e{assum3}, it follows then
\be \label{01}
|l(t)-l_e| = |(1- \f tT) (l^0-l_e) +\f tT (l^1-l_e)| \leq 2 \nu \leq 2 \nu_1.
\ee

\begin{figure}[htbp]
\includegraphics[width=0.45\textwidth]{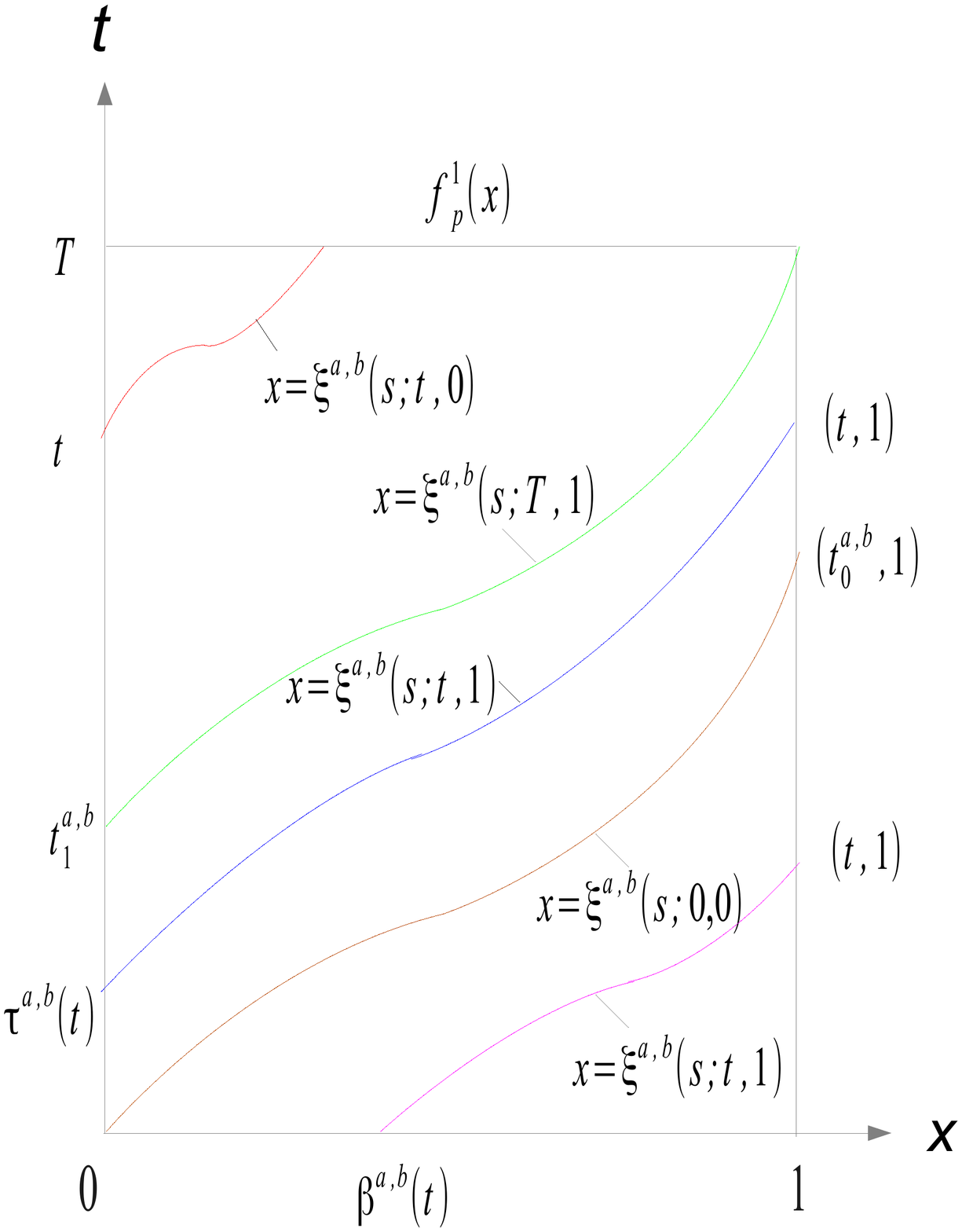}
\centering
\caption{The characteristics $\xi^{a,b}(s;t,x)$ and  $t_0^{a,b}$, $\beta^{a, b}(t)$, $\tau^{a, b}(t)$, $t_1^{a,b}$. \label{fig3} }\par\vspace{0pt}
\end{figure}

Next let us construct the desired control $F_{in}^{a,b}(t)$. For this purpose, we shall give the expression of $f_p(t,x)$ by characteristic method.
We denote by $\xi^{a,b}(s;t,x)$ the characteristic passing through $(t,x) \in [0,T] \times [0,1]$:
\be \label{xi-ab}
\begin{cases}
\displaystyle\f{d\xi^{a,b}(s;t,x)}{ds}= \alpha_p^{a,b} (s, \xi^{a,b}(s;t,x)) = \f{\zeta N^{a,b}(s) -  \xi^{a,b}(s;t,x) F^{a,b}(s)}{a(s)}, \\
\xi^{a,b}(t;t,x)=x.
\end{cases}
\ee
Suppose that the characteristic $\xi^{a,b}(s;0,0)$ intersects the line $x=1$ at $(t_0^{a,b},1)$.
Then for every $t\in[0,t_0^{a,b}]$,
the characteristic passing through $(t,1)$ intersects with $x$-axis at $(0,\beta^{a,b}(t))$;
for every $t\in(t_0^{a,b},T]$, the characteristic passing through $(t,1)$ intersects with $t$-axis at $(\tau^{a,b}(t),0)$.
The characteristic passing through $(T,1)$ intersects with $t$-axis at $(t_1^{a,b},0)$ (see Fig. \ref{fig3}).
Hence, we have
\begin{align*}
\xi^{a,b}(t^{a,b}_0;0,0)=1,\quad  \xi^{a,b}(0;t,1)=\beta^{a,b}(t),\quad
\xi^{a,b}(\tau^{a,b}(t);t,1)=0, \quad \xi^{a,b}(t_1^{a,b};T,1)=0.
\end{align*}

Solving \eqref{xi-ab} and using \eqref{talpha_p} and \eqref{FF}, we obtain
\begin{align}
\label{ua000}
 & e^{\int_{0}^{t^{a,b}_0} \f{l^1-l^0}{T a(\sigma)}\, d\sigma}
-\int_0^{t^{a,b}_0} \f{\zeta N^{a,b}(\sigma)}{a(\sigma)}\cdot e^{\int_0^{\sigma} \f{l^1-l^0}{T a(s)}\, ds}\, d\sigma =1,\\
\label{ua00}
 & e^{\int_0^t \f{l^1-l^0}{T a(\sigma)}\, d\sigma}
-\int_0^t \f{\zeta N^{a,b}(\sigma)}{a(\sigma)}\cdot e^{\int_0^{\sigma} \f{l^1-l^0}{T a(s)}\, ds}\, d\sigma =\beta^{a,b}(t),\\
\label{ua0}
& e^{\int_{\tau^{a,b}(t)}^t \f{l^1-l^0}{T a(\sigma)}\, d\sigma}
-\int_{\tau^{a,b}(t)}^t \f{\zeta N^{a,b}(\sigma)}{a(\sigma)}\cdot e^{\int_{\tau^{a,b}(t)}^{\sigma} \f{l^1-l^0}{T a(s)}\, ds}\, d\sigma =1,\\
\label{uaa}
& e^{\int_{t_1^{a,b}}^T \f{l^1-l^0}{T a(\sigma)}\, d\sigma}
-\int_{t_1^{a,b}}^T \f{\zeta N^{a,b}(\sigma)}{a(\sigma)}\cdot e^{\int_{t_1^{a,b}}^{\sigma} \f{l^1-l^0}{T a(s)}\, ds}\, d\sigma =1.
\end{align}
We define the control function $F^{a,b}_{in}(t)$ as the following
\be\label{cN20}
F^{a,b}_{in}(t)=\begin{cases}
h(t)\rho_0 V_{eff}N^{a,b}(t), \quad t\in [0,t_1^{a,b}],\\
f^1_p(\xi^{a,b}(T;t,0))\rho_0 V_{eff} N^{a,b}(t), \quad t\in [t_1^{a,b},T],
\end{cases}
\ee
where $h(t)$ is any artificial $W^{1,\infty}$ function satisfying
\be\label{hnu}
\|h(\cdot)-f_{pe}\|_{W^{1,\infty}}\leq \nu_1
\ee
and the following compatibility conditions
\begin{align}\label{h0}
h(0)= f^0_p(0),\quad
h(t_1^{a,b})= f^1_p(1).
\end{align}
Inspired by \eqref{sol1}, we define a map $\mathcal{F}:=(\mathcal{F}_1,\mathcal{F}_2): \Omega^{\varepsilon_1,T}\rightarrow C^0([0,T])$,
$\Phi\mapsto \mathcal{F}(\Phi)$ as
\begin{align}\label{F1}
\mathcal{F}_1(\Phi)(t):=& l(t)=l_0+\f{l^1-l^0}{T}t, \quad t\in [0,T],\\
\label{F2}
\mathcal{F}_2(\Phi)(t):=& \begin{cases}
f^0_p(\beta^{a,b}(t)),\, t\in[0,t^{a,b}_0),\\
h(\tau^{a,b}(t)),\, t\in (t^{a,b}_0,T],
\end{cases}
\end{align}
where $t^{a,b}_0$, $\beta^{a,b}(t)$ and $\tau^{a,b}(t)$ are defined by \eqref{ua000}, \eqref{ua00} and \eqref{ua0} respectively.

Now we prove that  $\mathcal{F}$ is a contraction mapping on $\Omega^{\varepsilon_1,T}$ provided that $\nu_1>0$ is small.
Obviously, the existence of the fixed point  implies the existence of the desired control to the original nonlinear controllability problem.

Let $\nu_1\leq  \varepsilon_1/2$. By \eqref{assum3}, \eqref{01}, \eqref{hnu}, \eqref{F1} and \eqref{F2},
 $\mathcal{F}(\Phi)\in \Omega^{\varepsilon_1,T}$ for all  $\Phi \in \Omega^{\varepsilon_1,T}$.

For any $\bar\Phi=(\bar a,\bar b)\in \Omega^{\varepsilon_1,T}$, we denote by $\xi^{\bar a,\bar b}(s;t,x)$ the characteristic
passing through $(t,x)\in [0,T]\times[0,1]$   as in \eqref{xi-ab} upon replacing $a,b$ by $\bar a, \bar b$. Correspondingly, we define
$t_0^{\bar a,\bar b}$, $\beta^{\bar a,\bar b}(t)$, $\tau^{\bar a,\bar b}(t)$, $t_1^{\bar a,\bar b}$
as  in \eqref{ua000}, \eqref{ua00}, \eqref {ua0} and \eqref{uaa} upon replacing $a,b$ by $\bar a, \bar b$.

By definition of $\mathcal{F}=(\mathcal{F}_1,\mathcal{F}_2)$ in \eqref{F1} and \eqref{F2},
we get  $\mathcal{F}_1(\bar\Phi) \equiv \mathcal{F}_1(\Phi) $ and thus
\begin{align*}
\|\mathcal{F}(\bar\Phi)-\mathcal{F}(\Phi)\|_{C^0([0,T])}
=\|\mathcal{F}_2(\bar\Phi)-\mathcal{F}_2(\Phi)\|_{C^0([0,T])}
=\sup_{t\in [0,T]}  |(\mathcal{F}_2(\bar\Phi)-\mathcal{F}_2(\Phi))(t)|.
\end{align*}

Without loss of generality, we may assume that $t^{\bar a,\bar b}_0>t^{a, b}_0$.
Hence, we need to estimate the point-wisely  $|\mathcal{F}_2(\bar\Phi)(t)-\mathcal{F}_2(\Phi)(t)|$
on the time interval $[0,t^{a,b}_0]$, $[t^{a,b}_0,t^{\bar a,\bar b}_0]$, $[t^{\bar a,\bar b}_0,T]$
respectively.

For any given $t\in[0,t^{a,b}_0]$, by \eqref{assum3} and \eqref{F2}, we have
\begin{align*} 
|(\mathcal{F}_2(\bar\Phi)-\mathcal{F}_2(\Phi))(t)|
=& \Big|f^0_p(\beta^{\bar a,\bar b}(t))-f^0_p(\beta^{a,b}(t))\Big|\nonumber\\
\leq &  \|f^0_{px}\|_{L^{\infty}} | \beta^{\bar a,\bar b}(t)- \beta^{a,b}(t)| \nonumber\\
\leq &  \nu_1 | \beta^{\bar a,\bar b}(t)- \beta^{a,b}(t)|.
\end{align*}
By \eqref{ua00}, it is easy to get
\begin{align*} 
| \beta^{\bar a,\bar b}(t)- \beta^{a,b}(t)|
\leq C (\|\bar a-a\|_{C^0([0,T])}+\|\bar b-b\|_{C^0([0,T])})
\leq C  \|\bar\Phi-\Phi\|_{C^0([0,T])}.
\end{align*}
Here and hereafter, we denote by $C$ various constants that are independent of $(\bar \Phi, \Phi)$.
Therefore,
\begin{align}\label{6.20}
\sup_{t\in [0, t^{a,b}_0]}|(\mathcal{F}_2(\bar \Phi)-\mathcal{F}_2(\Phi))(t)|
\leq C \nu_1  \|\bar\Phi-\Phi\|_{C^0([0,T])}.
\end{align}

For any given $t\in[t^{\bar a,\bar b}_0,T]$, by \eqref{hnu} and \eqref{F2}, we obtain
\begin{align}\label{es5}
|(\mathcal{F}_2(\bar \Phi)-\mathcal{F}_2(\Phi))(t)|
=|h(\tau^{\bar a,\bar b}(t))-h(\tau^{a,b}(t))|\leq \nu_1 |\tau^{\bar a,\bar b}(t)-\tau^{a,b}(t)|.
\end{align}
By \eqref{xi-ab}, for every $t\in[t^{a, b}_0,T]$,
\begin{align*}
\xi^{a,b}(s;t,1)=e^{\int_s^t \f{l^1-l^0}{T a(\sigma)}\, d\sigma}
-\int_s^t \f{\zeta N^{ a, b}(\sigma)}{ a(\sigma)}\cdot e^{\int_s^{\sigma} \f{l^1-l^0}{T  a(s)}\, ds}\, d\sigma, \quad s\in [\tau^{a,b}(t),T].
\end{align*}
Therefore,
\begin{align}  \label{es2}
  \sup_{s \in [\tau^{\bar a,\bar b}(t),T]} |\xi^{\bar a,\bar b}(s;t,1)-\xi^{a,b}(s;t,1)|
\leq C(\|\bar a-a\|_{C^0[0,T]}+\|\bar b-b\|_{C^0[0,T]}).
\end{align}
Then, for any $t\in [t^{\bar a,\bar b}_0, T]$, we have
\begin{align}\label{es3}
|\tau^{\bar a,\bar b}(t)-\tau^{a,b}(t)|
\leq &  \f{1}{\inf  \alpha_p^{a,b}} \left| \int_{\tau^{a,b}(t)}^{\tau^{\bar a,\bar b}(t)} \alpha_p ^{a,b}(s; \xi^{a,b}(s;t,1))\, ds \right|\nonumber\\ \nonumber
=& \f{1}{\inf  \alpha_p^{a,b}} \left| \int_{\tau^{\bar a,\bar b}(t)}^{T} \Big( \alpha_p ^{a,b}(s; \xi^{a,b}(s;t,1)) -
         \alpha_p ^{\bar a,\bar b}(s; \xi^{\bar a,\bar b}(s;t,1)) \Big) \, ds \right| \\ \nonumber
\leq &C\Big(\|\bar a-a\|_{C^0([0,T])}
+\|\bar b-b\|_{C^0([0,T])}\Big)
\\
\leq & C   \|\bar\Phi-\Phi\|_{C^0([0,T])},
\end{align}
where
\beq
\inf \alpha_p^{a,b} := \inf_{(t,x)\in [0,T]\times [0,1] \atop (a,b) \in \Omega^{\varepsilon_1,T} } \alpha_p^{a,b}(t,x) >0
\eeq
is independent of $\Phi$.
Therefore,
\begin{align}\label{6.25}
\sup_{t\in [t^{\bar a,\bar b}_0, T]}  |(\mathcal{F}_2(\bar \Phi)-\mathcal{F}_2(\Phi))(t)|
\leq C \nu_1  \|\bar\Phi-\Phi\|_{C^0([0,T])}.
\end{align}

For any given $t\in [t^{a,b}_0,t^{\bar a,\bar b}_0]$, by \eqref{h0} and \eqref{F2}, we have
\begin{align}\label{es02}
|\mathcal{F}_2(\bar\Phi)-\mathcal{F}_2(\Phi)|
=& |f^0_p(\beta^{\bar a,\bar b}(t))-h(\tau^{a,b}(t))|\nonumber\\
\leq & |f^0_p(\beta^{\bar a,\bar b}(t))-f^0_p(0)|+|h(\tau^{a,b}(t))-h(0)|\nonumber\\
\leq &\|f^0_{p_x}\|_{L^{\infty}}|\beta^{\bar a,\bar b}(t)|+\|h_t\|_{L^{\infty}}|\tau^{a,b}(t)| \nonumber\\
\leq  & \nu_1 (|\beta^{\bar a,\bar b}(t)|+|\tau^{a,b}(t)| ).
\end{align}
Note that   for any $ t\in [t^{a,b}_0,t^{\bar a,\bar b}_0]$,
\begin{align*}
\xi^{\bar a,\bar b}(s;t,1)=&e^{\int_s^t \f{l^1-l^0}{T \bar a(\sigma)}\, d\sigma}
-\int_s^t \f{\zeta N^{\bar a,\bar b}(\sigma)}{\bar a(\sigma)}\cdot e^{\int_s^{\sigma} \f{l^1-l^0}{T \bar a(s)}\, ds}\, d\sigma,
\quad s \in [0, t].
\end{align*}
Then, using $\xi^{\bar a,\bar b}(0;t^{\bar a,\bar b}_0,1)=0$, it follows that
\be\label{es04}
|\beta^{\bar a,\bar b}(t)|=|\xi^{\bar a,\bar b}(0;t,1)-\xi^{\bar a,\bar b}(0;t^{\bar a,\bar b}_0,1)|
\leq C|t-t^{\bar a,\bar b}_0|
\leq C|t^{\bar a,\bar b}_0-t^{a,b}_0|.
\ee
On the other hand, by \eqref{ua0} and note that $\tau^{a,b}(t_0^{a,b})=0$,
we can prove that for any $t\in[t^{a,b}_0,t^{\bar a,\bar b}_0]$,
\be\label{es05}
|\tau^{a,b}(t)|
= |\tau^{a,b}(t)-\tau^{a,b}(t_0^{a,b})|
\leq C  |t-t^{a,b}_0|
\leq C  |t^{\bar a,\bar b}_0-t^{a,b}_0|.
\ee

By definition of $t^{\bar a,\bar b}_0$ and $t^{a,b}_0$,
similarly to the derivation of \eqref{es2} and \eqref{es3}, we get
\begin{align}\label{es03}
|t^{\bar a,\bar b}_0-t^{ a, b}_0|\leq C(\|\bar a-a\|_{C^0([0,T])}+\|\bar b-b\|_{C^0([0,T])})
\leq  C   \|\bar\Phi-\Phi\|_{C^0([0,T])}.
\end{align}
Then it follows from  \eqref{es02}, \eqref{es04}, \eqref{es05} and \eqref{es03}, that
\begin{align} \label{6.30}
\sup_{t\in [t^{a,b}_0,t^{\bar a,\bar b}_0]} |(\mathcal{F}_2(\bar\Phi)-\mathcal{F}_2(\Phi))(t)|
\leq C\nu_1\|\bar\Phi-\Phi\|_{C^0([0,T])}.
\end{align}

Combining \eqref{6.20}, \eqref{6.25} and \eqref{6.30}, we can choose $\nu_1$ small enough
such that
\be
\|\mathcal{F}(\bar\Phi)-\mathcal{F}(\Phi)\|_{C^0([0,T])}
=\|\mathcal{F}_2(\bar\Phi)-\mathcal{F}_2(\Phi)\|_{C^0([0,T])}\leq \f{1}{2}\|\bar\Phi-\Phi\|_{C^0([0,T])}.
\ee
Then the contraction mapping theorem implies that $\mathcal{F}$ has a unique fixed point $(l(\cdot),f_p(\cdot,1))$ in $\Omega^{\varepsilon_1,T}$.
Therefore we find a solution $(l,f_p)$ to Cauchy problem \eqref{eq-filling-ratio-bound}-\eqref{eql-normal}
which also satisfies the final conditions  \eqref{final-con1}-\eqref{final-con2}.
Moreover, because of the   uniqueness of solution,
 the desired control function $N(t)$ and $F_{in}(t)$ can be chosen by substituting the fixed point $(l(\cdot), f_p(\cdot,1))$ into  \eqref{cN} and
\eqref{cN20} respectively. This concludes the proof of Theorem \ref{thm3}. \qed



\section{Acknowledgements}
The authors would like to thank Professor Jean-Michel Coron,  Professor Miroslav Krstic  and Professor Bernhard Maschke for their helpful comments and constant support.
The authors are thankful to  the support  of the ERC advanced grant 266907 (CPDENL) and the hospitality of the Laboratoire Jacques-Louis Lions of Universit\'{e} Pierre et Marie Curie.

Mamadou Diagne  has been supported by the French National Research Agency sponsored project ANR-11-BS03-0002 HAMECMOPSYS as a PhD candidate at  LAGEP (Laboratoire d'Automatique et du G\'enie des Proc\'ed\'es) of the university Claude Bernard Lyon I  and is currently supported by the  CCSD (Cymer Center for Control Systems and Dynamics) at the university of California San Diego as a postdoctoral fellow. Peipei Shang was partially supported by the National Science Foundation of China (No. 11301387) and by the Specialized Research Fund for the Doctoral Program of Higher Education (No. 20130072120008).
Zhiqiang Wang was partially supported by the National Science Foundation of China (No. 11271082) and by
the State Key Program of National Natural Science Foundation of China (No. 11331004).

\bibliographystyle{plain}


\begin{thebibliography}{10}

\bibitem{Stefano05}
S.~Bianchini and A.~Bressan.
\newblock Vanishing viscosity solutions of nonlinear hyperbolic systems.
\newblock {\em Ann. of Math.}, 161(1):223--342, 2005.

\bibitem{Bouchemal06}
K.~Bouchemal and  F.~Couenne and  S.~Briancon and  H.~Fessi, and M.~Tayakout.
\newblock Polyamides nanocapsules : Modeling and wall thickness estimation.
\newblock  {\em AIChE Journal}, 52:2161--2170, 2006.

\bibitem{BressanBook}
A.~ Bressan.
\newblock {\em Hyperbolic systems of conservation laws}, volume~20 of {\em
  Oxford Lecture Series in Mathematics and its Applications}.
\newblock Oxford University Press, Oxford, 2000.
\newblock The one-dimensional Cauchy problem.

\bibitem{CoronBook}
J.-M. Coron.
\newblock {\em Control and nonlinearity}, volume 136 of {\em Mathematical
  Surveys and Monographs}.
\newblock American Mathematical Society, Providence, RI, 2007.

\bibitem{coron}
J.-M. Coron.
\newblock Local controllability of a 1-{D} tank containing a fluid
              modeled by the shallow water equations.
\newblock {\em ESAIM Control Optim. Calc. Var.}, 8:513--554, 2002.




\bibitem{CGWang}
J.-M. Coron and O. Glass, and Z. Wang.
\newblock Exact boundary controllability for 1-{D} quasilinear
              hyperbolic systems with a vanishing characteristic speed.
\newblock {\em SIAM J. Control Optim.}, 48:3150--3122, 2010.

\bibitem{CKW}
J.-M. Coron and  M. Kawski, and Z. Wang.
\newblock Analysis of a conservation law modeling a highly re-entrant
  manufacturing system.
\newblock {\em Discrete Contin. Dyn. Syst. Ser. B}, 14(4):1337--1359, 2010.




\bibitem{CW2012}
J.-M. Coron and Z.  Wang.
\newblock Controllability for a scalar conservation law with nonlocal velocity.
\newblock {\em J. Differential Equations}, 252:181--201, 2012.

\bibitem{Gugat}
M. Gugat and G. Leugering.
\newblock Global boundary controllability of the de {S}t.\ {V}enant
              equations between steady states.
\newblock {\em Ann. Inst. H. Poincar\'e Anal. Non Lin\'eaire.}, 20:1--11, 2003.

\bibitem{Daraoui2010}
N.~Daraoui and P.~Dufour and H.~Hammouri, and A.~Hottot.
\newblock Model predictive control during the primary drying stage of
  lyophilisation.
\newblock {\em Control Engineering Practice}, 18(5):483--494, 2010.

\bibitem{Diagne}
M. Diagne and  V.~Dos Santos~Martins and  F.~Couenne and B.~Maschke.
\newblock Well posedness of the model of an extruder in infinite dimension.
\newblock In {\em Decision and Control and European Control Conference
  (CDC-ECC), 2011 50th IEEE Conference on}, pages 1311--1316, 2011.

\bibitem{DiagneJESA11}
M. Diagne and V.~Dos Santos~Martins and  F.~Couenne and  B.~Maschke, and C.~Jallut.
\newblock Mod{\'e}lisation et commande d'un syst{\`e}me d'{\'e}quations aux
  d{\'e}riv{\'e}es partielles \`a fronti{\`e}re mobile : application au
  proc{\'e}d{\'e} d'extrusion.
\newblock {\em Journal Europ{\'e}en des Syst{\`e}mes Automatis{\'e}s},
  45:665--691, 2011.


\bibitem{DSW2015}
M. Diagne and P. Shang  and Z. Wang.
\newblock Feedback stabilization for the mass balance equations of an extrusion
  process.
\newblock {\em IEEE Transactions on Automatic Control}, 2015, doi={10.1109/TAC.2015.2444232}.

\bibitem{DSW2015b}
M. Diagne and P. Shang  and Z. Wang.
\newblock Feedback Stabilization of a Food Extrusion
Process Described by 1D PDEs Defined on
Coupled Time-Varying Spatial Domains
\newblock {\em 12th IFAC Workshop on Time Delay Systems}, June 28-30, 2015 Ann Arbor, MI, USA, pages 51--56, 2015.

\bibitem{Godlewski_ESAIM05}
E. Godlewski and  K.-C. Le~Than and P.A. Raviart.
\newblock The numerical interface coupling of nonlinear systems of conservation
  laws: Ii. the case of systems.
\newblock {\em ESAIM: Mathematical Modelling and Numerical Analysis},
  39(4):649--692, 2005.

\bibitem{Godlewski04}
E. Godlewski and P.-A. Raviart.
\newblock The numerical interface coupling of nonlinear hyperbolic systems of
  conservation laws: I. the scalar case.
\newblock {\em Numer. Math.}, 97:81--130, 2004.

\bibitem{Horma3}
L.  H{\"o}rmander.
\newblock {\em The analysis of linear partial differential operators. {III}},
  volume 274 of {\em Grundlehren der Mathematischen Wissenschaften [Fundamental
  Principles of Mathematical Sciences]}.
\newblock Springer-Verlag, Berlin, 1994.
\newblock Pseudo-differential operators, Corrected reprint of the 1985
  original.

\bibitem{KIM001}
E.~K. Kim and J.~L. White.
\newblock Isothermal transient startup for starved flow modular co-rotating
  twin screw extruder.
\newblock {\em Polymer Engineering and Science}, 40:543--553, 2004.

\bibitem{KIM002}
E.~K. Kim and J.~L. White.
\newblock Non-isothermal transient startup for starved flow modular co-rotating
  twin screw extruder.
\newblock {\em International Polymer Processing}, 15:233--241, 2004.

\bibitem{KULSH92}
M. K. Kulshrestha and C.A. Zaror.
\newblock An unsteady state model for twin screw extruders.
\newblock {\em Tran IChemE, PartC}, 70:21--28, 1992.

\bibitem{LeFlochBook}
P.~G. LeFloch.
\newblock {\em Hyperbolic systems of conservation laws}.
\newblock Lectures in Mathematics ETH Z\"urich. Birkh\"auser Verlag, Basel,
  2002.
\newblock The theory of classical and nonclassical shock waves.

\bibitem{CHIN01}
C.~H.  Li.
\newblock Modelling extrusion cooking.
\newblock {\em Mathematical and Computer Modelling}, 33:553--563, 2001.



\bibitem{LiBook09}
T.-T. Li.
\newblock {\em Controllability and observability for quasilinear hyperbolic
              systems}.
\newblock American Institute of Mathematical Sciences (AIMS),
              Springfield, MO, 2010.

\bibitem{LiBook94}
T.-T. Li.
\newblock {\em Global classical solutions for quasilinear hyperbolic systems}.
\newblock Research in Applied Mathematics {\bf 32}, John Wiley \& Sons,
  Chichester, 1994.


\bibitem{Li2001}
T.-T.  Li and Y.~Jin.
\newblock Semi-global $c_1$ solution to the mixed initial-boundary value
  problem for quasilinear hyperbolic systems.
\newblock {\em Chinese Ann. Math. Ser. B}, 22:325--336, 2001.

\bibitem{LiRao}
T.-T.  Li and B. Rao.
\newblock Exact boundary controllability for quasi-linear hyperbolic
              systems.
\newblock {\em SIAM J. Control Optim.}, 41:1748--1755, 2003.

\bibitem{LiYuBook}
T.-T.  Li and W. Yu.
\newblock {\em Boundary value problems for quasilinear hyperbolic systems}.
\newblock Duke University Mathematics Series, V. Duke University Mathematics
  Department, Durham, NC, 1985.

\bibitem{LiuYang}
T. Liu and T. Yang.
\newblock Well-posedness theory for hyperbolic conservation laws.
\newblock {\em Comm. Pure Appl. Math.}, 52(12):1553--1586, 1999.

\bibitem{Petit10}
N.~Petit.
\newblock Control problems for one-dimensional fluids and reactive fluids with
  moving interfaces.
\newblock In {\em Advances in the theory of control, signals and systems with
  physical modeling}, volume 407 of {\em Lecture notes in control and
  information sciences}, pages 323--337, Lausanne, Dec 2010.

\bibitem{Purlis10}
E.~Purlis and V.~O. Salvadori.
\newblock A moving boundary problem in a food material undergoing volume change-Simulation of bread baking.
\newblock {\em Food Research International}, 43:949--958, 2008.

\bibitem{Russell}
D. L. Russell.
\newblock Controllability and stabilizability theory for linear partial
              differential equations: recent progress and open questions.
\newblock {\em SIAM Rev.}, 20:639--739, 1978.

\bibitem{SW}
P. Shang and Z.  Wang.
\newblock Analysis and control of a scalar conservation law modeling a highly
  re-entrant manufacturing system.
\newblock {\em J. Differential Equations}, 250(2):949--982, 2011.

\bibitem{Velardi08}
S.~A. Velardi and A.~A. Barresi.
\newblock Development of simplified models for the freeze-drying process and
  investigation of the optimal operating conditions.
\newblock {\em Chemical Engineering Research and Design}, 86:9--22, 2008.

\bibitem{W2006}
Z.  Wang.
\newblock Exact controllability for nonautonomous first order quasilinear
  hyperbolic systems.
\newblock {\em Chinese Ann. Math. Ser. B}, 27:643--656, 2006.



\end{thebibliography}

\end{document}